\documentclass[a4paper,11pt]{article}

\input{diagrams}

\newcommand{\op}{\ensuremath{\mbox{\hspace{1pt}{\scriptsize op}}}}
\newcommand{\set}{\ensuremath{\mbox{\bfseries Set}}}
\newcommand{\Set}{\ensuremath{\mbox{\bfseries Set}}}
\newcommand{\cat}[1]{\ensuremath{\mbox{\bfseries {\upshape {#1}}}}}
\newcommand{\free}[1]{\ensuremath{{\mathcal F} \hspace{-1.5pt}
    {#1}^{\op}}}

\newcommand{\elt}[1]{\ensuremath{\mbox{\upshape elt} \hspace{1pt} {#1}}}

\newcommand{\cl}[1]{\ensuremath{\mathcal {#1}}}
\newcommand{\bb}[1]{\ensuremath{\mathbb {#1}}}

\newcommand{\lra}{\ensuremath{\longrightarrow}}
\newcommand{\map}[1]{\ensuremath{\stackrel{{#1}}{\lra}}}
\newcommand{\simarr}{\ensuremath{\stackrel{\sim}{\longrightarrow}}}
\newcommand{\sm}{symmetric multicategory{}}
\newcommand{\sms}{symmetric multicategories{}}

\newcommand{\iso}{isomorphism{}}

\newcommand{\oset}{\cat{OSet}}

\newcommand{\pshfcat}[1]{\ensuremath{[{#1}^{\op},\cat{Set}]}}

\newtheorem{theorem}{Theorem}[section]
\newtheorem{proposition}[theorem]{Proposition}
\newtheorem{corollary}[theorem]{Corollary}
\newtheorem{lemma}[theorem]{Lemma}
\newtheorem{definition}[theorem]{Definition}
\newtheorem{definitions}[theorem]{Definitions}

\newenvironment{prf}{\vspace{2ex}\begin{sloppypar}{\noindent\upshape
{\bfseries Proof. }}} {{\hspace*{\fill}
$\Box$}\end{sloppypar}\vspace{2ex}}

\newenvironment{prfof}[1]{\vspace{2ex}\begin{sloppypar}{\noindent
\upshape{\bfseries Proof of {#1}. }}} {{\hspace*{\fill}
$\Box$}\end{sloppypar}\vspace{2ex}}

\newcommand{\numarabic}{\renewcommand{\labelenumi}{\arabic{enumi})}}

\newcommand{\pica}{\begin{center} \input}
\newcommand{\picz}{\end{center}}

\newlength{\leng}
\newlength{\fontleng}
\newcommand{\sunit}{\setlength{\unitlength}{1mm}}
\newcommand{\setleng}[1]{\setlength{\leng}{#1}
    \setlength{\unitlength}{\leng}}
\newcommand{\setunit}[1]{\setlength{\unitlength}{#1}}

\sunit

\sunit

\sunit

\sunit

\sunit

\sunit

\sunit

\sunit

\sunit

\newcommand{\scalecq}{

\put(0,0){
\setlength{\unitlength}{0.3\leng} %

\put(20,90){\line(0,1){20}}      %
\put(40,90){\line(0,1){20}}      %
\put(80,90){\line(0,1){20}}      %

\put(20,90){\line(1,0){60}}      %
\put(20,90){\line(1,-1){30}}     %
\put(80,90){\line(-1,-1){30}}    %
\put(50,40){\line(0,1){20}}      %
\put(60,104){\makebox[0pt]{$\cdots$}} }}

\newcommand{\scalecr}[5]{

\put(0,0){
\setunit{0.3\leng} %
\scalecq
\put(20,114){\makebox[0pt]{#1}}   %
\put(40,114){\makebox[0pt]{#2}}   %
\put(80,114){\makebox[0pt]{#3}}   %
\put(50,76){\makebox(0,0){#5}}      
\put(45,38){\makebox(0,0)[t]{#4}}  }} 

\usepackage{amsfonts}
\usepackage{amssymb}
\usepackage{amsmath}
\usepackage{fancyheadings}
\usepackage{latexcad}

\begin{document}

\title{The category of opetopes and the category of opetopic sets}

\author{Eugenia Cheng\\ \\Department of Pure Mathematics, University
of Cambridge\\E-mail: e.cheng@dpmms.cam.ac.uk}
\date{October 2002}
\maketitle

\begin{abstract} 

We give an explicit construction of the category
\cat{opetope} of opetopes.  We prove that the category of opetopic
sets is equivalent to the category of presheaves over \cat{opetope}

\end{abstract}

\setcounter{tocdepth}{3}
\tableofcontents

\section*{Introduction}
\addcontentsline{toc}{section}{Introduction}

In \cite{bd1}, Baez and Dolan give a definition of weak $n$-category
in which the underlying shapes of cells are `opetopes' and the
underlying data is given by `opetopic sets'.  The idea is that
opetopic sets should be presheaves over the category of opetopes. 
However Baez and Dolan do not explicitly construct the category of
opetopes, so opetopic sets are defined directly instead.  A
relationship between this category of opetopic sets and a category of
presheaves is alluded to but not proved.  

The main result of this paper is that the category of opetopic sets is
equivalent to the category of presheaves over the category of
opetopes.  However, we do not use the opetopic definitions exactly as
given in \cite{bd1} but continue to use the modifications given in our
earlier work (\cite{che7}, \cite{che8}).  In these papers we use a
generalisation along lines which the original authors began, but chose
to abandon for reasons which are unclear.  This generalisation enables
us, in \cite{che7}, to exhibit a relationship with the work of Hermida,
Makkai and Power (\cite{hmp1}) and, in \cite{che8}, with the work of
Leinster (\cite{Lei5}).  Given these useful results, we continue to
study the modified theory in this work.

We begin, in Section~\ref{opecat} by giving an explicit construction
of the category of opetopes.  The idea is as follows.  In \cite{che7} we constructed, for each $k
\geq 0$, a category $\bb{C}_k$ of $k$-opetopes.  For the category
$\cat{Opetope}$ of opetopes of all dimensions, each category
$\bb{C}_k$ should be a full subcategory of \cat{Opetope}; furthermore
there should be `face maps' exhibiting the constituent $m$-opetopes,
or `faces' of a $k$-opetope, for $m \le k$. We refer to the
$m$-opetope faces as $m$-faces.  Note that there are no degeneracy maps.  

The $(k-1)$-faces of a $k$-opetope $\alpha$ should be the
$(k-1)$-opetopes of its source and target; these should all be
distinct.  Then each of these faces has its own $(k-2)$-faces, but
all these $(k-2)$-opetopes should not necessarily be considered as
distinct $(k-2)$-faces in $\alpha$. For $\alpha$ is a
configuration for composing its $(k-1)$-faces at their
$(k-2)$-faces, so the $(k-2)$-faces should be identified with one
another at places where composition is to occur.  That is, the
composite face maps from these $(k-2)$-opetopes to $\alpha$ should
therefore be equal.  Some further details are then required to
deal with isomorphic copies of opetopes.

Recall that a `configuration' for composing $(k-1)$-opetopes is
expressed as a tree (see \cite{che7}) whose nodes
are labelled by the $(k-1)$-opetopes in question, with the edges
giving their inputs and outputs.  So composition occurs along each
edge of the tree, via an object-morphism label, and thus the tree
tells us which $(k-1)$-opetopes are identified.

In order to express this more precisely, we first give a more
formal description of trees (Section~\ref{treeformal}).  In fact,
this leads to an abstract description of trees as certain
Kelly-Mac~Lane graphs; this is the subject of \cite{che12}.  The
results of Section~\ref{treeformal} thus arise as preliminary
results in \cite{che12} and we refer the reader to this paper for
the full account and proofs.  

In Section~\ref{osets}, we examine the theory of opetopic sets. 
We begin by following through our modifications to the opetopic
theory to include the theory of opetopic sets.  (Our previous work
has only dealt with the theory of opetopes.)  We then use results of
\cite{kel3} to prove that the category of opetopic sets is indeed
equivalent to the category of presheaves on $\cl{O}$, the category of
opetopes defined in Section~\ref{opecat}.  This is the main result of
this work.  

Finally, a comment is due on the notion of `multitope' as
defined in \cite{hmp1}.  In this work, Hermida, Makkai and Power
begin a definition of $n$-category explicitly analogous to that of
\cite{bd1}, the analogous concepts being `multitopes' and `multitopic
sets'.  In \cite{che7} we prove that `opetopes and multitopes are the
same up to isomorphism', that is, for each $k \geq 0$ the category of
$k$-opetopes is equivalent to the (discrete) category of
$k$-multitopes.  In \cite{hmp1}, Hermida, Makkai and Power do go on
to give an explicit definition of the analogous category
\cat{Multitope}, of multitopes.  Given the above equivialences, and
assuming the underlying idea is the same, this would be equivalent to
the category \cat{Opetope}, but we do not attempt to prove it in this
work.

\subsubsection*{Terminology}

\renewcommand{\labelenumi}{\roman{enumi})}
\begin{enumerate}

\item Since we are concerned chiefly with {\em weak}
$n$-categories, we follow Baez and Dolan (\cite{bd1}) and omit the
word `weak' unless emphasis is required; we refer to {\em strict}
$n$-categories as `strict $n$-categories'.

\item We use the term `weak $n$-functor' for an $n$-functor where
functoriality holds up to coherent isomorphisms, and `lax' functor
when the constraints are not necessarily invertible.

\item In \cite{bd1} Baez and Dolan use the terms `operad' and
`types' where we use `multicategory' and `objects'; the latter
terminology is more consistent with Leinster's use of `operad' to
describe a multicategory whose `objects-object' is 1.

\item In \cite{hmp1} Hermida, Makkai and Power use the term
`multitope' for the objects constructed in analogy with the
`opetopes' of \cite{bd1}.  This is intended to reflect the fact
that opetopes are constructed using operads but multitopes using
multicategories, a distinction that we have removed by using the
term `multicategory' in both cases.  However, we continue to use
the term `opetope' and furthermore, use it in general to refer to
the analogous objects constructed in each of the three theories.
Note also that Leinster uses the term `opetope' to describe
objects which are analogous but not {\em a priori} the same; we
refer to these as `Leinster opetopes' if clarification is needed.

\item We regard sets as sets or discrete categories with no
notational distinction.

\end{enumerate}

\bigskip {\bfseries Acknowledgements}

This work was supported by a PhD grant from EPSRC.  I would like to thank
Martin Hyland and Tom Leinster for their support and guidance.

\section{The category of opetopes}
\label{opecat}

In this section we give an explicit construction of the category
\cat{Opetope} of opetopes.  This construction will enable us, in
Section~\ref{osets}, to prove that the category of opetopic sets
is in fact a presheaf category.

We begin with a brief account of the trees used to construct
higher-dimensional opetopes from lower-dimensional ones; we refer the
reader to \cite{che12} for the full account, with proofs and examples.

\sunit
\subsection{Informal description of trees} \label{treeback}

Recall the trees introduced in \cite{che7} to
describe the morphisms of a slice multicategory.  These are
`labelled combed trees' with ordered nodes.  In fact, we will
first consider the {\em unlabelled} version of such trees, since
the labelled version follows easily.  For example the following is
a tree:

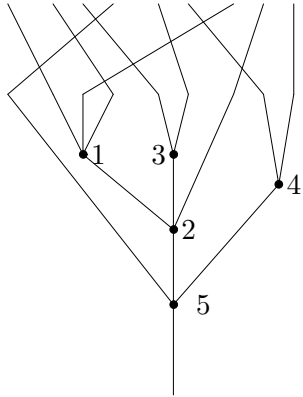
\begin{center}
\begin{picture}(48,60)
\thinlines
\drawthickdot{26.0}{16.0}
\drawthickdot{26.0}{26.0}
\drawthickdot{26.0}{36.0}
\drawthickdot{14.0}{36.0}
\drawthickdot{40.0}{32.0}
\drawpath{14.0}{36.0}{10.0}{44.0}
\drawpath{14.0}{36.0}{14.0}{44.0}
\drawpath{14.0}{36.0}{18.0}{44.0}
\drawpath{26.0}{36.0}{24.0}{44.0}
\drawpath{26.0}{36.0}{28.0}{44.0}
\drawpath{26.0}{26.0}{14.0}{36.0}
\drawpath{26.0}{26.0}{26.0}{36.0}
\drawpath{26.0}{26.0}{34.0}{44.0}
\drawpath{26.0}{16.0}{26.0}{26.0}
\drawpath{26.0}{16.0}{26.0}{4.0}
\drawpath{26.0}{16.0}{4.0}{44.0}
\drawpath{40.0}{32.0}{38.0}{44.0}
\drawpath{40.0}{32.0}{42.0}{44.0}
\drawpath{26.0}{16.0}{40.0}{32.0}
\drawpath{10.0}{44.0}{4.0}{56.0}
\drawpath{18.0}{44.0}{10.0}{56.0}
\drawpath{14.0}{56.0}{24.0}{44.0}
\drawpath{18.0}{56.0}{4.0}{44.0}
\drawpath{24.0}{56.0}{28.0}{44.0}
\drawpath{28.0}{56.0}{38.0}{44.0}
\drawpath{34.0}{56.0}{14.0}{44.0}
\drawpath{38.0}{56.0}{34.0}{44.0}
\drawpath{42.0}{56.0}{42.0}{44.0}
\drawcenteredtext{16.0}{36.0}{1}
\drawcenteredtext{28.0}{26.0}{2}
\drawcenteredtext{24.0}{36.0}{3}
\drawcenteredtext{42.0}{32.0}{4}
\drawcenteredtext{30.0}{16.0}{5}
\end{picture}

\end{center}

\renewcommand{\labelenumi}{\roman{enumi})}

Explicitly, a tree $T=(T, \rho, \tau)$ consists of
\begin{enumerate}
\item A planar tree $T$
\item A permutation $\rho \in \cat{S}_l$ where $l=$ number of leaves of
$T$
\item A bijection $\tau: \{\mbox{nodes of $T$}\} \lra \{1, 2,
\ldots, k\}$ where $k=$ number of nodes of $T$; equivalently an
ordering on the nodes of $T$.
\end{enumerate}


Note that there is a `null tree' with no nodes
\begin{center}
\begin{picture}(14,24)
\thinlines
\drawpath{4.0}{20.0}{4.0}{4.0}
\end{picture}
 .
\end{center}

%

\subsection{Formal description of trees}
\label{treeformal}

In this section we give a formal description of the above trees,
characterising them as connected graphs with no closed loops (in
the conventional sense of `graph').  This will enable us, in
Section~\ref{loopy}, to determine which faces of faces are
identified in an opetope.

Note that the material in this section is presented fully in
\cite{che12}.  It enables us to express a tree as a Kelly-Mac~Lane
graph; it also enables us to show that all
{\em allowable} Kelly-Mac~Lane graphs of the correct shape arise in this way.

We consider a tree with $k$ nodes $N_1, \ldots, N_k$ where $N_i$
has $m_i$ inputs and one output. Let $N$ be a node with
$(\sum\limits_i m_i )- k+1$ inputs; $N$ will be used to represent
the leaves and root of the tree.

Then a tree is given by a bijection
    \[\coprod_i\{\mbox{inputs of $N_i$}\} \coprod \{\mbox{output of $N$}\} \lra
    \coprod_i\{\mbox{output of $N_i$}\} \coprod \{\mbox{inputs of
    $N$}\}\]
since each input of a node is either connected to a unique output
of another node, or it is a leaf, that is, input of $N$. Similarly
each output of a node is either attached to an input of another
node, or it is the root, that is, output of $N$.

We express this formally as follows.

\begin{lemma} \label{tree1}
Let $T$ be a tree with nodes $N_1, \ldots, N_k$, where $N_i$ has
inputs $\{x_{i1}, \ldots, x_{im_i}\}$ and output $x_i$.  Let $N$
be a node with inputs $\{z_1, \ldots, z_l\}$ and output $z$, with
    \[l=(\sum\limits_{i=1}^k m_i) - k +1.\]
    Then $T$ is given by a
bijection
    \[\alpha \ : \ \coprod_i\{x_{i1},\ldots, x_{im_i}\} \coprod \{z\} \lra
    \coprod_i\{x_i\} \coprod \{z_1, \ldots, z_l\}.\]
\end{lemma}

For the converse, every such bijection gives a graph, but it is not
necessarily a tree.  We need to
ensure that the resulting graph has no closed loops; the use of
the `formal' node $N$ then ensures connectedness.  We express this
formally as follows.

\begin{lemma} \label{loopspropb}
Let $N_1, \ldots, N_k, N$ be nodes where $N_i$ has inputs
$\{x_{i1}, \ldots, x_{im_i}\}$ and output $x_i$, and $N$ has
inputs $\{z_1, \ldots, z_l\}$ and output $z$, with
$l=(\sum\limits_{i=1}^k m_i) - k +1$. Let $\alpha$ be a bijection
    \[\coprod_i\{x_{i1},\ldots, x_{im_i}\} \coprod \{z\} \lra
    \coprod_i\{x_i\} \coprod \{z_1, \ldots, z_l\}.\]
Then $\alpha$ defines a graph with nodes $N_1, \ldots, N_k$.
\end{lemma}

\begin{lemma} \label{loopspropb2}
Let $\alpha$ be a graph as above.  Then $\alpha$ has a closed loop
if and only if there is a non-empty sequence of indices
    \[\{t_1 , \ldots, t_n \} \subseteq \{1, \ldots, k\}\]
such that for each $2 \leq j \leq n$
    \[\alpha( x_{t_jb_j})=x_{t_{j-1}} \]
for some $1 \leq b_j \leq m_j$, and
    \[\alpha(x_{t_1b_1})= x_{t_n}\]
for some $1 \leq b_1 \leq m_1$.  \end{lemma}

\begin{corollary} \label{loopspropb3}
A tree with nodes $N_1, \ldots, N_k$ is precisely a bijection
$\alpha$ as in Lemma~\ref{loopspropb}, such that there is no
sequence of indices as in Lemma~\ref{loopspropb2}.
\end{corollary}

\subsection{Labelled trees}
\label{labeltree}

For the construction of opetopes we require the `labelled' version
of the trees presented in Section~\ref{treeback}.  A tree labelled
in a category \bb{C} is a tree as above, with each edge labelled
by a morphism of \bb{C} considered to be pointing `down' towards
the root.

\begin{proposition} Let $N_1, \ldots, N_k, N$ be nodes where $N_i$ has inputs
\[\{x_{i1}, \ldots, x_{im_i}\}\] and output $x_i$, and $N$ has inputs
$\{z_1, \ldots, z_l\}$ and output $z$, with
\[l=(\sum\limits_{i=1}^k m_i) - k +1.\] Then a labelled tree with
these nodes is given by a bijection
    \[\alpha: \coprod_i\{x_{i1},\ldots, x_{im_i}\} \coprod \{z\} \lra
    \coprod_i\{x_i\} \coprod \{z_1, \ldots, z_l\}\]
satisfying the conditions as above, together with, for each
    \[y \in  \coprod_i\{x_{i1},\ldots, x_{im_i}\} \coprod \{z\}\]
a morphism $f \in \bb{C}$ giving the label of the edge joining $y$
and $\alpha(y)$.  Then $y$ is considered to be labelled by the
object $\mbox{{\upshape cod}}(f)$ and $\alpha(y)$ by the object
$\mbox{{\upshape dom}}(f)$.
\end{proposition}

\begin{prf} Follows immediately from Corollary~\ref{loopspropb3} and
the definition. \end{prf}

\subsection{The category of opetopes} \label{loopy}

In our earlier work (\cite{che7}) we constructed for each $k\geq 0$ the
category $\bb{C}_k$ of $k$-opetopes.  We now construct a category
\cat{Opetope} of opetopes of all dimensions whose morphisms are,
essentially, face maps.  Each category $\bb{C}_k$ is to be a full
subcategory of \cat{Opetope}, and there are no morphisms from an
opetope to one of lower dimension.

We construct the category $\cat{Opetope} = \cl{O}$ as follows.
Write $\cl{O}_k = \bb{C}_k$.

For the objects:
    \[\mbox{ob\ }\cl{O} = \coprod_{k\geq 0} \cl{O}_k.\]

The morphisms of \cl{O} are given by generators and relations as
follows.

\begin{itemize}
\item Generators
\end{itemize}

\numarabic
\begin{enumerate}

\item For each morphism $f: \alpha \lra \beta \in \cl{O}_k$ there
is a morphism
    \[f: \alpha \lra \beta \in \cl{O}.\]

\item Let $k \geq 1$ and consider $\alpha \in \cl{O}_k = o(I^{k+}) =
\elt{(I^{(k-1)+})}$.  Write $\alpha \in I^{(k-1)+}(x_1, \ldots,
x_m;x)$.  Then for each $1 \leq i \leq m$ there is a morphism
    \[s_i: x_i \lra \alpha \in \cl{O}\]
and there is also a morphism
    \[t: x \lra \alpha \in \cl{O}.\]
We write $G_k$ for the set of all generating morphisms of this
kind.

\end{enumerate}

Before giving the relations on these morphisms we make the
following observation about morphisms in $\cl{O}_k$.  Consider
    \begin{eqnarray*} \alpha & \in & I^{(k-1)+}(x_1, \ldots, x_m; x)\\
    \beta & \in & I^{(k-1)+}(y_1, \ldots, y_m; y) \end{eqnarray*}
A morphism $\alpha \map{g} \beta \in \cl{O}_k$ is given by a
permutation $\sigma$ and morphisms
    \[\begin{array}{ccccl} x_i & \map{f_i} & y_{\sigma(i)} &&\\
    x & \map{f} & y & \in & \cl{O}_{k-1} \end{array}\]
So for each face map $\gamma$ there is a unique `restriction' of
$g$ to the specified face, giving a morphism $\gamma g$ of
$(k-1)$-opetopes.

Note that, to specify a morphism in the category
$\free{\cl{O}_{k-1}} \times \cl{O}_{k-1}$ the morphisms $f_i$
above should be in the direction $y_{\sigma(i)} \lra x_i$, but
since these are all unique isomorphisms the direction does not
matter; the convention above helps the notation.  We now give the
relations on the above generating morphisms.

\begin{itemize}
\item{Relations}
\end{itemize}

\begin{enumerate}

\item For any morphism
    \[\alpha \map{g} \beta  \ \in \cl{O}_k\]
and face map
    \[x_i \map{s_i} \alpha\]
the following diagrams commute
    \[\begin{diagram} x_i & \rTo^{s_i} & \alpha && x & \rTo^{t} & \alpha \\
    \dTo{s_i(g)}{} & & \dTo{}{g} &&  \dTo{t(g)}{} & & \dTo{}{g} \\
    y_{\sigma(i)} & \rTo_{s_{\sigma(i)}} & \beta && y & \rTo_{t} & \beta \end{diagram}\]

We write these generally as
    \[\begin{diagram} x & \rTo^{\gamma} & \alpha \\
    \dTo{\gamma g}{} & & \dTo{}{g} \\
    y & \rTo_{\gamma'} & \beta \end{diagram}\]

\item Faces are identified where composition occurs: consider
$\theta \in \cl{O}_k$ where $k\geq 2$.  Recall that
$\theta$ is constructed as an arrow of a slice multicategory, so
is given by a labelled tree, with nodes labelled by its
$(k-1)$-faces, and edges labelled by object-morphisms, that is,
morphisms of $\cl{O}_{k-2}$.

So by the formal description of trees (Section~\ref{treeformal}),
$\theta$ is a certain bijection, and the elements that are in
bijection with each other are the $(k-2)$-faces of the
$(k-1)$-faces of $\theta$; they are given by composable pairs of
face maps of the second kind above.  That is, the node labels are
given by face maps $\alpha \map{\gamma} \theta$ and then the
inputs and outputs of those are given by pairs
    \[x \map{\gamma_1} \alpha \map{\gamma_2} \theta\]
where $\gamma_2 \in G_k$ and $\gamma_1 \in G_{k-1}$.  Now, if
\[\begin{array}{rccccc}  & x & \map{\gamma_1} & \alpha &
    \map{\gamma_2} & \theta \\
    \mbox{and } & y & \map{\gamma_3} & \beta & \map{\gamma_4} & \theta
\end{array}\]
correspond under the bijection, there must be a unique
object-morphism
    \[f:x \lra y\]
labelling the relevant edge of the tree.  Then for the composites
in $\cl{O}$ we have the relation: the following diagram commutes

    \[\begin{diagram}
        x     & \rTo^{\gamma_1} &   \alpha   &             & \\
    \dTo{f}{} &      &            & \rdTo(2,1)^{\gamma_2}  & \theta \\
        y     & \rTo^{\gamma_3} &    \beta   & \ruTo(2,1)_{\gamma_4}  & \end{diagram}.\]

\item Composition in $\cl{O}_k$ is respected, that is, if $g \circ
f = h \in \cl{O}_k$ then $g \circ f = h \in \cl{O}$.

\item Identities in $\cl{O}_k$ are respected, that is, given any
morphism $x \map{\gamma} \alpha \in \cl{O}$ we have $\gamma \circ
1_x = \gamma$.

\end{enumerate}

Note that only the relation (2) is concerned with the
identification of faces with one another; the other relations are
merely dealing with isomorphic copies of opetopes.

We immediately check that the above relations have not identified
any morphisms of $\cl{O}_k$.

\begin{lemma} Each $\cl{O}_k$ is a full subcategory of \cl{O}.
\end{lemma}

\begin{prf}  Clear from definitions.  \end{prf}

We now check that the above relations have not identified any
$(k-1)$-faces of $k$-opetopes.

\begin{proposition}\label{normal} Let $x \in \cl{O}_{k-1}$, $\alpha \in
\cl{O}_k$ and $\gamma_1, \gamma_2 \in G_k$ with
    \[\gamma_1,\ \gamma_2 : x \lra \alpha\]
Then $\gamma_1 = \gamma_2 \in \cl{O} \implies \gamma_1 = \gamma_2
\in G_k$.  \end{proposition}

We prove this by expressing all morphisms from $(k-1)$-opetopes to
$k$-opetopes in the following ``normal form''; this is a simple
exercise in term rewriting (see \cite{kv1}).

\begin{lemma} Let $x \in \cl{O}_{k-1}$, $\alpha \in \cl{O}$.  Then
a morphism
    \[x \lra \alpha \ \in \cl{O}\]
is uniquely represented by
    \[x \map{\gamma} \alpha\]
or a pair
    \[x \map{f} y \map{\gamma} \alpha\]
where $f \in \cl{O}_{k-1}$ and $\gamma \in G_k$. \end{lemma}

\begin{prf}  Any map $x \lra \alpha$ is represented by terms of
the form
    \[x \map{f_1} x_1 \map{f_2} \cdots \map{f_m} x_m \map{\gamma}
    \alpha_1 \map{g_1} \cdots \map{g_{j-1}}
    \alpha_j \map{g_j} \alpha\]
where each $f_i \in \cl{O}_{k-1}$ and each $g_r \in \cl{O}_k$.
Equalities are generated by equalities in components of the
following forms: \numarabic
\[\begin{diagram}
    1) &{}& \rTo^{\gamma} &{}& \rTo^{g} & = & \rTo^{\gamma g} &{}&
    \rTo^{\gamma'}&{} \\
    2) &{}& \rTo^f &{}& \rTo^{f'} & = & \rTo^{f' \circ f} &{} & \in
    \cl{O}_{k-1} & \\
    3) &{}& \rTo^g &{}& \rTo^{g'} & = & \rTo^{g' \circ g} &{} & \in
    \cl{O}_{k} & \\
    4) &{}& \rTo^1 &{}& \rTo^{\gamma} & = & \rTo^{\gamma} &{} & &
\end{diagram}\]
where $\gamma \in G_k$ and $\gamma g$ and $\gamma'$ are as defined
above.  That is, equalities in terms are generated by equations $t
= t'$ where $t'$ is obtained from $t$ by replacing a component of
$t$ of a left hand form above, with the form in the right hand
side, or vice versa.

We now orient the equations in the term rewriting style in the
direction
    \[\begin{diagram} {}& \implies &{} \end{diagram}\]
from left to right in the above equations.  We then show two
obvious properties:

\begin{enumerate}

\item Any reduction of $t$ by $\Longrightarrow$ terminates in at
most $2j+m$ steps.

\item If we have \pica pic76.tex \picz then there exists $t'''$ with
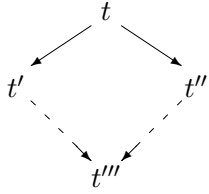
\pica pic75.tex \picz where the dotted arrows indicate a chain of
equations (in this case of length at most 2).
\end{enumerate}

The first part is clear from the definitions; for the second part
the only non-trivial case is for a component of the form
    \[\begin{diagram} {}& \rTo^{\gamma} &{}& \rTo^{g_1} &{}&
    \rTo^{g_2} &{} \end{diagram}.\]
This reduces uniquely to
    \[\begin{diagram} {}& \rTo^{\gamma(g_2 \circ g_1)} & {} & \rTo
    ^{\gamma'} &{} \end{diagram}\]
since `restriction' is unique, as discussed earlier.

It follows that, for any terms $t$ and $s$, $t=s$ if and only if
$t$ and $s$ reduce to the same normal form as above.  \end{prf}

\begin{prfof}{Proposition \ref{normal}  } $\gamma_1$ and $\gamma_2$
are in normal form.  \end{prfof}


\section{Opetopic Sets} \label{osets}

In this section we examine the theory of opetopic sets. We begin by
following through our modifications to the opetopic theory to include
the theory of opetopic sets.  We then use results of \cite{kel3} to
prove that the category of opetopic sets is indeed equivalent to the
category of presheaves on $\cl{O}$, the category of opetopes defined
in Section~\ref{opecat}.

Recall that, by the equivalences proved in the \cite{che7} and
\cite{che8}, we have equivalent categories of opetopes, multitopes
and Leinster opetopes.  So we may define equivalent categories of
opetopic sets by taking presheaves on any of these three categories. 
In the following definitions, although the opetopes we consider are
the `symmetric multicategory' kind, the concrete description of an
opetopic set is not {\em precisely} as a presheaf on the category of
these opetopes.  The sets given in the data are indexed not by
opetopes themselves but by {\em isomorphism classes} of opetopes; so
at first sight this resembles a presheaf on the category of Leinster
opetopes.  However, we do not pursue this matter here, since the
equivalences proved in our earlier work are sufficient for the
purposes of this article.

We adopt this presentation in order to avoid naming the same cells
repeatedly according to the symmetries; that is, we do not keep
copies of cells that are isomorphic by the symmetries.

\subsection{Definitions}
\label{osetnewdef}

In \cite{bd1}, weak $n$-categories are defined as opetopic sets
satisfying certain universality conditions. However, opetopic sets
are defined using only \sms\ with a {\em set} of objects; in the
light of the results of our earlier work, we seek a
definition using \sms\ with a {\em category} of objects. The
definitions we give here are those given in \cite{bd1} but with
modifications as demanded by the results of our previous work.

\label{bicatopicsets}

The underlying data for an opetopic $n$-category are given by an
opetopic set. Recall that, in \cite{bd1}, given a symmetric multicategory $Q$ 
a $Q$-opetopic set $X$ is given by, for each $k \geq 0$, a \sm\ $Q(k)$ and a 
set $X(k)$ over $o(Q(k))$, where
\[\begin{array}{c}
    Q(0) = Q\\
    \mbox{and \ \ } Q(k+1) = {Q(k)_{X(k)}}^+.
\end{array}\]  An opetopic set is then an $I$-opetopic set, where $I$ is the
\sm\ with one object and one (identity) arrow.

The idea is that the category of opetopic sets should be
equivalent to the presheaf category \[[\mbox{{\bfseries
Opetope}}^{\op},\mbox{{\bfseries Set}}]\] and we use this to
motivate our generalisation of the Baez-Dolan definitions.



Recall that we have for each $k \geq 0$ a category $\bb{C}(k)$
of $k$-opetopes, and each $\bb{C}(k)$ is a full subcategory of
{\bf Opetope}. A functor \[\mbox{{\bfseries Opetope}}^{\op}
\longrightarrow \mbox{{\bfseries Set}}\] may be considered as
assigning to each opetope a set of `labels'.

Recall that for each $k$, $\bb{C}(k)$ is equivalent to a discrete
category.  So it is sufficient to specify `labels' for each
isomorphism class of opetopes. 

Recall (\cite{che7}) that we call a \sm\ $Q$ {\em tidy} if
it is freely symmetric with a category of objects \bb{C}
equivalent to a discrete category.  Throughout this section we say
`$Q$ has object-category \bb{C} equivalent to $S$ discrete' to
mean that $S$ is the set of isomorphism classes of \bb{C}, so
\bb{C} is equipped with a morphism $\bb{C} \simarr S$.  We begin
by defining the construction used for `labelling' as discussed
above.  The idea is to give a set of labels as a set over the
isomorphism classes of objects of $Q$, and then to `attach' the
labels using the following pullback construction.

\begin{definition} Let $Q$ be a tidy symmetric multicategory with
category of objects \bb{C} equivalent to $S$ discrete.  Given a
set $X$ over $S$, that is, equipped with a function $f:X \lra S$,
we define the {\em pullback multicategory} $Q_X$ as follows.

\begin{itemize}

\item Objects: $o(Q_X)$ is given by the pullback

\setlength{\unitlength}{0.2em}
\begin{center}
\begin{picture}(45,35)(5,7)       %

\put(10,10){\makebox(0,0){$X$}}  
\put(10,35){\makebox(0,0){.}}  
\put(45,35){\makebox(0,0){$\bb{C}$}}  
\put(45,10){\makebox(0,0){$S$}}  

\put(15,35){\vector(1,0){27}}  
\put(15,10){\vector(1,0){27}}  
\put(10,32){\vector(0,-1){19}} 
\put(45,32){\vector(0,-1){19}} 

\put(29,8){\makebox(0,0)[t]{$f$}}  %
\put(50,22){\makebox(0,0){$\sim$}} %
\put(52,9){\makebox(0,0){.}}

\end{picture}
\end{center}  Observe that the morphism on the left is an
equivalence, so $o(Q_X)$ is equivalent to $X$ discrete.  Write $h$
for this morphism.

\item Arrows: given objects $a_1, \ldots a_k, a \in o(Q_X)$ we
have
    \[Q_X(a_1, \ldots, a_k; a) \cong Q(fh(a_1), \ldots, fh(a_k);
    fh(a)).\]

\item Composition, identities and symmetric action are then
inherited from $Q$.

\end{itemize}
\end{definition}

We observe immediately that since $Q$ is tidy, $Q_X$ is tidy. Also
note that if $Q$ is object-discrete this definition corresponds to
the definition of pullback symmetric multicategory given in
\cite{bd1}.

We are now ready to describe the construction of opetopic sets.

\begin{definition}  Let $Q$ be a tidy \sm\ with object-category
\bb{C} equivalent to $S$ discrete.  A $Q$-{\em opetopic set} $X$
is defined recursively as a set $X(0)$ over $S$ together with a
${Q_X}^+$-opetopic set $X_1$.
\end{definition}

So a $Q$-opetopic set consists of, for each $k \geq 0$:
\begin{itemize}
\item a tidy \sm\ $Q(k)$ with object-category $\bb{C}(k)$
equivalent to $S(k)$ discrete
\item a set $X(k)$ and function $X(k)
\stackrel{f_k}{\longrightarrow} S(k)$
\end{itemize}
where
\begin{eqnarray*} Q(0) &=& Q\\
     \mbox{and \ \ } Q(k+1) &=& {Q(k)_{X(k)}}^+. \end{eqnarray*} We refer to
$X_1$ as the underlying ${Q(k)_{X(k)}}^+$-opetopic set of $X$.

We now define morphisms of opetopic sets.  Suppose we have
opetopic sets $X$ and $X'$ with notation as above, together with a
morphism of symmetric multicategories
    \[F: Q \lra Q'\]
and a function
    \[F_0: X(0) \lra X'(0)\]
such that the following diagram commutes

\setlength{\unitlength}{0.2em}
\begin{center}
\begin{picture}(45,40)(5,5)        %

\put(10,10){\makebox(0,0){$X'(0)$}}  
\put(10,35){\makebox(0,0){$X(0)$}}  
\put(45,35){\makebox(0,0){$S(0)$}}  
\put(45,10){\makebox(0,0){$S'(0)$}}  

\put(17,35){\vector(1,0){21}}  
\put(17,10){\vector(1,0){21}}  
\put(10,30){\vector(0,-1){15}} 
\put(45,30){\vector(0,-1){15}} 

\put(8,23){\makebox(0,0)[r]{$F_0$}} 
\put(47,23){\makebox(0,0)[l]{$F$}} 
\put(27,37){\makebox(0,0)[b]{$f_0$}} 
\put(27,8){\makebox(0,0)[t]{$f'_0$}} 


\end{picture}
\end{center}where the morphism on the right is given by the action of $F$ on
objects. This induces a morphism
    \[Q_{{X(0)}} \lra Q'_{X'(0)}\]
and so a morphism
    \[{Q_{{X(0)}}}^+ \lra {Q'_{{X'(0)}}}^+.\]
We make the following definition.

\begin{definition} A {\em morphism of $Q$-opetopic sets}
    \[F: X \lra X'\]
is given by:

\begin{itemize}

\item an underlying morphism of \sms\ and function $F_0$ as above

\item a morphism $X_1 \lra X'_1$ of their underlying opetopic
sets, whose underlying morphism is induced as above.

\end{itemize}
\end{definition}

\sunit

\noindent So $F$ consists of
\begin{itemize}
\item a morphism $Q \lra Q'$
\item for each $k \geq 0$ a function $F_k : X(k) \lra X'(k)$ such
that the following diagram commutes
\end{itemize}

\setlength{\unitlength}{0.2em}

\begin{center}
\begin{picture}(45,40)(5,5)      %

\put(10,10){\makebox(0,0){$X'(k)$}}  
\put(10,35){\makebox(0,0){$X(k)$}}  
\put(45,35){\makebox(0,0){$S(k)$}}  
\put(45,10){\makebox(0,0){$S'(k)$}}  

\put(17,35){\vector(1,0){21}}  
\put(17,10){\vector(1,0){21}}  
\put(10,30){\vector(0,-1){15}} 
\put(45,30){\vector(0,-1){15}} 

\put(8,23){\makebox(0,0)[r]{$F_k$}} 
\put(47,23){\makebox(0,0)[l]{}} 
\put(27,37){\makebox(0,0)[b]{$f_k$}} 
\put(27,8){\makebox(0,0)[t]{$f'_k$}} 


\end{picture}
\end{center} where the map on the right hand side is induced as appropriate.

Note that the above notation for a $Q$-opetopic set $X$ and
morphism $F$ will be used throughout this section, unless
otherwise specified.

\begin{definition} An {\em opetopic set} is an $I$-opetopic set.
A morphism of opetopic sets is a morphism of $I$-opetopic sets. We
write \cat{OSet} for the category of opetopic sets and their
morphisms.
\end{definition}

Eventually, a weak $n$-category is defined as an opetopic set with
certain properties.  The idea is that $k$-cells have underlying
shapes given by the objects of $I^{k+}$.  These are `unlabelled'
cells.  To make these into fully labelled $k$-cells, we first give
labels to the 0-cells, via the function $X(0) \longrightarrow
S(0)$, and then to 1-cells via $X(1) \longrightarrow S(1)$, and so
on.  This idea may be captured in the following `schematic'
diagram.

\sunit
\begin{picture}(140,100)(0,-10)
\put(0,20){\makebox(0,0)[c]{${Q(2)_{{X(2)}}}^+$}}
\put(0,40){\makebox(0,0)[c]{$Q(2)_{{X(2)}}$}}
\put(37,40){\makebox(0,0)[c]{$Q(2) = {Q(1)_{{X(1)}}}^+$}}
\put(43,60){\makebox(0,0)[c]{${Q(1)_{{X(1)}}}$}}
\put(77,60){\makebox(0,0)[c]{$Q(1)={Q(0)_{{X(0)}}}^+$}}
\put(83,80){\makebox(0,0)[c]{$Q(0)_{{X(0)}}$}}
\put(113,80){\makebox(0,0)[c]{$Q(0)=I$}}
\put(120,60){\makebox(0,0)[c]{$I^+$}}
\put(120,40){\makebox(0,0)[c]{$I^{2+}$}}
\put(120,20){\makebox(0,0)[c]{$I^{3+}$}}

\put(119,76){\vector(0,-1){12}}  
\put(119,56){\vector(0,-1){12}}  %
\put(119,36){\vector(0,-1){12}}  %
\put(119,16){\vector(0,-1){12}}  %
\put(79,76){\vector(0,-1){12}}  %
\put(39,56){\vector(0,-1){12}}  %
\put(-1,36){\vector(0,-1){12}}  %

\put(121,70){\makebox(0,0)[l]{$+$}}  %
\put(121,50){\makebox(0,0)[l]{$+$}}  %
\put(121,30){\makebox(0,0)[l]{$+$}}  %
\put(121,10){\makebox(0,0)[l]{$+$}}  %
\put(81,70){\makebox(0,0)[l]{$+$}}  %
\put(41,50){\makebox(0,0)[l]{$+$}}  %
\put(1,30){\makebox(0,0)[l]{$+$}}  %

\put(92,80){\vector(1,0){11}}   
\put(51,60){\vector(1,0){10}}   %
\put(9,40){\vector(1,0){11}}   %

\put(60,0){$\vdots$}

\end{picture}


Bearing in mind our modified definitions, we use the Baez-Dolan
terminology as follows.

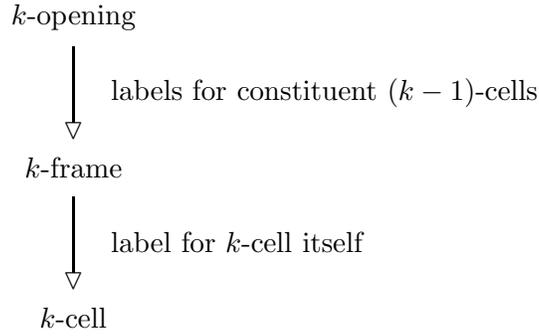
\begin{definitions} \
\begin{itemize}
\item A $k$-{\em dimensional cell} (or $k$-{\em cell}) is an element
of $X(k)$ \\ (i.e. an isomorphism class of objects of
$Q(k)_{{X(k)}}$\,).
\item A $k$-{\em frame} is an isomorphism class of objects of
$Q(k)$ \\ (i.e. an \iso\ class of arrows of
$Q(k-1)_{{X(k-1)}}$\,).
\item A $k$-{\em opening} is an isomorphism class of arrows of
$Q(k-1)$, for $k \geq 1$.
\end{itemize}
\end{definitions}

So a $k$-opening may acquire $(k-1)$-cell labels and become a
$k$-frame, which may itself acquire a label and become a $k$-cell.
We refer to such a cell and frame as being {\em in} the original
$k$-opening.

On objects, the above schematic diagram becomes:

\sunit

\begin{picture}(140,120)(0,-30)
\put(15,20){\makebox(0,0)[c]{3-frames}}
\put(15,40){\makebox(0,0)[c]{2-cells}}
\put(50,20){\makebox(0,0)[c]{3-openings}}  %
\put(50,40){\makebox(0,0)[c]{2-frames}}  %
\put(50,60){\makebox(0,0)[c]{1-cells}}  %
\put(84,20){\makebox(0,0)[c]{.}}  %
\put(85,40){\makebox(0,0)[c]{2-openings}}  %
\put(85,60){\makebox(0,0)[c]{1-frames}}  %
\put(85,80){\makebox(0,0)[c]{0-cells}}  %
\put(120,20){\makebox(0,0)[c]{3-opetopes}}  %
\put(120,40){\makebox(0,0)[c]{2-opetopes}}  %
\put(120,60){\makebox(0,0)[c]{1-opetopes}}  %
\put(120,80){\makebox(0,0)[c]{0-opetopes}}  %

\put(119,76){\vector(0,-1){12}}  
\put(119,56){\vector(0,-1){12}}  %
\put(119,36){\vector(0,-1){12}}  %
\put(119,16){\vector(0,-1){12}}  %
\put(84,76){\vector(0,-1){12}}  %
\put(49,56){\vector(0,-1){12}}  %
\put(14,36){\vector(0,-1){12}}  %

\put(121,70){\makebox(0,0)[l]{$+$}}  %
\put(121,50){\makebox(0,0)[l]{$+$}}  %
\put(121,30){\makebox(0,0)[l]{$+$}}  %
\put(121,10){\makebox(0,0)[l]{$+$}}  %
\put(86,70){\makebox(0,0)[l]{$+$}}  %
\put(51,50){\makebox(0,0)[l]{$+$}}  %
\put(16,30){\makebox(0,0)[l]{$+$}}  %
\put(86,50){\makebox(0,0)[l]{$+$}}  %
\put(51,30){\makebox(0,0)[l]{$+$}}  %
\put(16,10){\makebox(0,0)[l]{$+$}}  %

\qbezier[15](84,56)(84,50)(84,45) 
\put(84,46){\vector(0,-1){2}}  %
\qbezier[15](84,36)(84,30)(84,25) %
\put(84,26){\vector(0,-1){2}}  %
\qbezier[15](49,36)(49,30)(49,25) %
\put(49,26){\vector(0,-1){2}}  %
\qbezier[15](84,16)(84,10)(84,5) %
\put(84,6){\vector(0,-1){2}}  %
\qbezier[15](49,16)(49,10)(49,5) %
\put(49,6){\vector(0,-1){2}}  %
\qbezier[15](14,16)(14,10)(14,5) %
\put(14,6){\vector(0,-1){2}}  %

\put(96,80){\vector(1,0){12}}   
\put(61,60){\vector(1,0){12}}   %
\put(26,40){\vector(1,0){12}}   %
\put(-9,20){\vector(1,0){12}}   %

\qbezier[15](96,60)(102,60)(107,60) 
\put(106,60){\vector(1,0){2}}  %
\qbezier[15](96,40)(102,40)(107,40) 
\put(106,40){\vector(1,0){2}}  %
\qbezier[15](96,20)(102,20)(107,20) 
\put(106,20){\vector(1,0){2}}  %
\qbezier[15](61,40)(67,40)(72,40) 
\put(71,40){\vector(1,0){2}}  %
\qbezier[15](61,20)(67,20)(72,20) 
\put(71,20){\vector(1,0){2}}  %
\qbezier[15](26,20)(32,20)(37,20) 
\put(36,20){\vector(1,0){2}}  %


\put(20,-10){\line(1,0){25}}   
\put(45,-10){\makebox(0,0)[l]{$\rhd$}}  %
\put(55,-10){\line(1,0){25}}   
\put(80,-10){\makebox(0,0)[l]{$\rhd$}}  %
\put(90,-10){\line(1,0){25}}   
\put(115,-10){\makebox(0,0)[l]{$\rhd$}}  %

\put(32.5,-15){\makebox(0,0)[c]{{\em labels for 2-cells}}}  %
\put(67.5,-15){\makebox(0,0)[c]{{\em labels for 1-cells}}}  %
\put(102.5,-15){\makebox(0,0)[c]{{\em labels for 0-cells}}}  %

\end{picture}.

Horizontal arrows represent the process of labelling, as shown;
vertical arrows represent the process of `moving up' dimensions.
Starting with a $k$-opetope, we have from right to left the
progressive labelling of 0-cells, 1-cells, and so on, to form a
$k$-cell at the far left, the final stages being:

\sunit
\begin{center}
\begin{picture}(70,60)(0,-10)

\put(10,40){\makebox(0,0)[c]{$k$-opening}}  %
\put(10,20){\makebox(0,0)[c]{$k$-frame}}  %
\put(10,0){\makebox(0,0)[c]{$k$-cell}}  %

\put(15,30){\makebox(0,0)[l]{labels for constituent $(k-1)$-cells}}  %
\put(15,10){\makebox(0,0)[l]{label for $k$-cell itself}}  %

\put(10,6){\line(0,1){10}}  %
\put(10,6){\makebox(0,0)[t]{$\triangledown$}}  %
\put(10,26){\line(0,1){10}}  %
\put(10,26){\makebox(0,0)[t]{$\triangledown$}}  %

\end{picture}
\end{center} A $k$-opening acquires labels as an arrow of $Q(k-1)$, becoming a
$k$-frame as an arrow of $Q(k-1)_{{X(k-1)}}$\,.  That is, it has
$(k-1)$-cells as its source and a $(k-1)$-cell as its target.

\begin{definition}
A $k$-{\em niche} is a $k$-opening (i.e. arrow of $Q(k-1)$)
together with labels for its source only.
\end{definition}

We may represent these notions as follows.  Let $f$ be an arrow of
$Q(k-1)$, so $f$ specifies a $k$-opening which we might represent
as

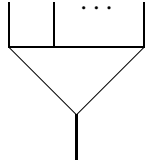
\setleng{1mm}
\begin{center}
\begin{picture}(18,30)(8,8)
\scalecq
\end{picture}.
\end{center}

\noindent Then a niche in $f$ is represented by

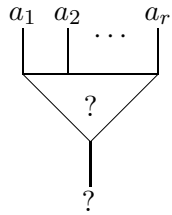
\begin{center}
\begin{picture}(18,30)(8,8)
\scalecr{$a_1$}{$a_2$}{$a_r$}{\ \ ?}{?}
\end{picture}
\end{center}

\noindent where $a_1, \ldots a_r$ are `valid' labels for the
source elements of $f$; a $k$-frame is represented by

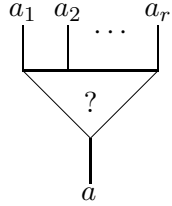
\begin{center}
\begin{picture}(18,30)(8,8)
\scalecr{$a_1$}{$a_2$}{$a_r$}{\ \ $a$}{?}
\end{picture}
\end{center}

\noindent where $a$ is a `valid' label for the target of $f$.
Finally a $k$-cell is represented by

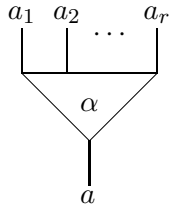
\begin{center}
\begin{picture}(18,30)(8,8)
\scalecr{$a_1$}{$a_2$}{$a_r$}{\ \ $a$}{$\alpha$}
\end{picture}.
\end{center}

Since all \sms\ in question are tidy, we may in each case
represent the same isomorphism class by any symmetric variant of
the above diagrams.  Also, we refer to $k$-cells as labelling
$k$-opetopes, rather than isomorphism classes of $k$-opetopes.

\subsection{\cat{OSet} is a presheaf category} \label{presheaf}

In this section we prove the main result of this work, that the
category of opetopic sets is a presheaf category, and moreover, that
it is equivalent to the presheaf category \[\pshfcat{\cl{O}}.\] To
prove this we use \cite{kel3}, Theorem 5.26, in the case $\cl{V} =
\Set$.  This theorem is as follows.

\begin{theorem} \label{kel} Let \cl{C} be a \cl{V}-category.
In order that $\cl{C}$ be equivalent to $[\cl{E}^{\op}, \cl{V}]$
for some small category $\cl{E}$ it is necessary and sufficient
that $\cl{C}$ be cocomplete, and that there be a set of
small-projective objects in $\cl{C}$ constituting a strong
generator for $\cl{C}$. \end{theorem}

We see from the proof of this theorem that if $E$ is such a set
and $\cl{E}$ is the full subcategory of $\cl{C}$ whose objects are
the elements of $E$, then
    \[\cl{C} \simeq [\cl{E}^{\op}, \cl{V}].\]
We prove the following propositions; the idea is to ``realise''
each isomorphism class of opetopes as an opetopic set; the set of
these opetopic sets constitutes a strong generator as required.

\begin{proposition} \label{pshf1}
\oset\ is cocomplete.
\end{proposition}

\begin{proposition} \label{pshf2}
There is a full and faithful functor
    \[G:\cl{O} \lra \oset.\]
\end{proposition}

\begin{proposition} \label{pshf3}
Let $\alpha \in \cl{O}$.  Then $G(\alpha)$ is small-projective in
\oset.
\end{proposition}

\begin{proposition} \label{pshf4}
Let
    \[E = \coprod\{G(\alpha)\ |\ \alpha \in \cl{O}\} \subseteq
    \oset.\]
Then $E$ is a strongly generating set for \oset.
\end{proposition}

\begin{corollary} \label{pshf5}
\oset\ is a presheaf category.
\end{corollary}

\begin{corollary} \label{pshf6}
    \[\oset \simeq \pshfcat{\cl{O}}.\]
\end{corollary}

\begin{prfof}{Proposition \ref{pshf1}}  Consider a diagram
    \[D: \bb{I} \lra \oset\]
where \bb{I} is a small category.  We seek to construct a limit
$Z$ for $D$; the set of cells of $Z$ of shape $\alpha$ is given by
a colimit of the sets of cells of shape $\alpha$ in each $D(I)$.

We construct an opetopic set $Z$ as follows.  For each $k \geq 0$,
$Z(k)$ is a colimit in \set:
    \[Z(k) = \int\nolimits^{I\in \bb{I}} D(I)(k).\]
Now for each $k$ we need to give a function
    \[F(k):Z(k) \lra o(Q(k))\]
where
    \[Q(k) = {Q(k-1)_{Z(k-1)}}^+\]
    \[Q(0)=I.\]
That is, for each $\alpha \in Z(k)$ we need to give its frame.
Now
    \[ \begin{array} {r} Z(k) = \coprod\limits_{I\in\bb{I}} D(I)(k) \\ \  \end{array}
    \left/ \begin{array}{l} \\ \sim \end{array} \right.  \]
where $\sim$ is the equivalence relation generated by
    \[ \begin{array} {cc} D(u)(\alpha_{I'})\sim \alpha_I & \mbox{\
    for all } u:I \lra I' \in \bb{I} \\
    & \mbox{and } \alpha_I \in D(I)(k). \end{array}\]
So $\alpha \in Z(k)$ is of the form $[\alpha_I]$ for some
$\alpha_I \in D(I)(k)$ where $[\alpha_I]$ denotes the equivalence
class of $\alpha_I$ with respect to $\sim$.

Now suppose the frame of $\alpha_I$ in $D(I)$ is
    \[(\beta_1, \ldots, \beta_j) \map{?} \beta\]
where $\beta_i, \beta \in D(I)(k-1)$ label some $k$-opetope $x$.
We set the frame of $[\alpha_I]$ to be
    \[(\ [\beta_1], \ldots, [\beta_j]\ ) \map{?} [\beta]\]
labelling the same opetope $x$.  This is well-defined since a
morphism of opetopic sets preserves frames of cells, so the frame
of $D(u)(\alpha_I)$ is
    \[(\ D(u)(\beta_1)\ , \ldots,\ D(u)(\beta_j)\ ) \map{?} D(u)(\beta)\]
also labelling $k$-opetope $x$.  It follows from the universal
properties of the colimits in \set\ that $Z$ is a colimit for $D$,
with coprojections induced from those in \set.  Then, since \set\
is cocomplete, \oset\ is cocomplete.  \end{prfof}

\begin{prfof}{Proposition \ref{pshf2}}  Let $\alpha$ be a
$k$-opetope.  We express $\alpha$ as an opetopic set
$G(\alpha)=\hat{\alpha}$ as follows, using the usual notation for
an opetopic set.  The idea is that the $m$-cells are given by the
$m$-faces of $\alpha$.

For each $m \geq 0$ set
    \[\begin{array}{cl} X(m) = \{\ [(x,f)]\ \  | & \mbox{ $x \in
    \cl{O}_m$ \ and $x \map{f} \alpha \in \cl{O}$}\\
    & \mbox{ where $[\ ]$ denotes isomorphism class in
    $\cl{O}/\alpha$} \}. \end{array}\]
So in particular we have
    \[X(k) = \{[(\alpha, 1)]\}\]
and for all $m>k$, $X(m) = \emptyset$. It remains to specify the
frame of $[(x,f)]$.  The frame is an object of
    \[Q(m) = {Q(m-1)_{X(m-1)}}^+\]
so an arrow of
    \[{Q(m-2)_{X(m-2)}}^+\]
labelled with elements of $X(m-1)$. Now such an arrow is a
configuration for composing arrows of $Q(m-2)_{X(m-2)}$; for the
frame as above, this is given by the opetope $x$ as a labelled
tree. Then the $(m-1)$-cell labels are given as follows. Write
    \[x :y_1, \ldots, y_j \lra y\]
say, and so we have for each $i$ a morphism
    \[y_i \lra x\]
and a morphism
    \[y \lra x \in \cl{O}.\]
Then the labels in $X(m-1)$ are given by
    \[[y_i \lra x \map{f} \alpha] \in X(m-1)\]
and
    \[[y \lra x \map{f} \alpha] \in X(m-1).\]
Now, given a morphism
    \[h: \alpha \lra \beta \in \cl{O}\]
we define
    \[\hat{h}: \hat{\alpha} \lra \hat{\beta} \in \cat{OSet}\]
by
    \[[(x,f)] \mapsto [(x, h \circ f)]\]
which is well-defined since if $(x,f) \cong (x', f')$ then $(x,
hf) \cong (x', hf')$ in $\cl{O}/\alpha$.  This is clearly a
morphism of opetopic sets.

Observe that any morphism $\hat{\alpha} \lra \hat{\beta}$ must be
of this form since the faces of $\alpha$ must be preserved.
Moreover, if $\hat{h} = \hat{g}$ then certainly
$[(\alpha,h)]=[(\alpha,g)]$.  But this gives $(\alpha,h) =
(\alpha,g)$ since there is a unique morphism $\alpha \lra \alpha
\in \cl{O}$ namely the identity. So $G$ is full and faithful as
required.  \end{prfof}

\begin{prfof}{Proposition \ref{pshf3}}  For any $\alpha \in \cl{O}_k$
we show that $\hat{\alpha}$ is small-projective, that is that the
functor
    \[\psi=\oset(\hat{\alpha},-): \oset \lra \set\]
preserves small colimits.  First observe that for any opetopic set
$X$
    \[ \begin{array}{ccl} \psi(X) = \oset(\hat{\alpha},X)& \cong &
    \{ \mbox{$k$-cells in $X$ whose underlying $k$-opetope is $\alpha$} \}\\
    & \subseteq & X(k) \end{array}\]
and the action on a morphism $F:X \lra Y$ is given by
    \[\begin{array}{cccc} \psi(F) = \oset(\hat{\alpha},F): & \oset(\hat{\alpha},X)&
    \lra & \oset(\hat{\alpha},Y) \\
    &x & \mapsto & F(x). \end{array}\]
So $\psi$ is the `restriction' to the set of cells of shape
$\alpha$.  This clearly preserves colimits since the cells of
shape $\alpha$ in the colimit are given by a colimit of the sets
cells of shape $\alpha$ in the original diagram. \end{prfof}

\begin{prfof}{Proposition \ref{pshf4}}  First note that
    \[\hat{\alpha} = \hat{\beta} \iff \alpha \cong \beta \in
    \cl{O}\]
so
    \[E \cong \coprod_k S_k\]
where for each $k$, $S_k$ is the set of $k$-dimensional Leinster
opetopes.  Since each $S_k$ is a set it follows that $E$ is a set.

We need to show that, given a morphism of opetopic sets $F:X \lra
Y$, we have
    \[\oset(\hat{\alpha},F) \mbox{ is an isomorphism for all
    $\hat{\alpha}$\ } \Longrightarrow \mbox{\ $F$ is an
    isomorphism.}\]
Now, we have seen above that
    \[\oset(\hat{\alpha}, X) \cong \{ \mbox{cells of $X$ of shape
     $\alpha$} \}\]
so
    \[\oset(\hat{\alpha},F) = F|_\alpha = \mbox{ $F$ restricted to
    cells of shape $\alpha$}.\]
So
    \[ \begin{array}{l} \oset(\hat{\alpha},F) \mbox{ is an isomorphism for all
    $\hat{\alpha}$\ \hspace{20mm}} \\
    \mbox{\hspace{20mm}} \iff  \mbox{ $F|_\alpha$ is an isomorphism for
    all $\alpha \in \cl{O}$} \\
    \mbox{\hspace{20mm}} \iff \mbox{$F$ is an isomorphism}. \end{array}\]
\end{prfof}

\begin{prfof}{Corollary \ref{pshf5}}  Follows from Propositions
\ref{pshf1}, \ref{pshf2}, \ref{pshf3}, \ref{pshf4} and \cite{kel3}
Theorem 5.26.  \end{prfof}

\begin{prfof}{Corollary \ref{pshf6}}  Let \cl{E} be the full
subcategory of \oset\ whose objects are those of $E$.  Since $G$
is full and faithful, $\cl{E}$ is the image of $G$ and we have
    \[\cl{O} \simeq \cl{E}\]
and hence
    \[\oset \simeq \pshfcat{\cl{E}} \simeq \pshfcat{\cl{O}}.\]
\end{prfof}

\addcontentsline{toc}{section}{References}
\bibliography{bib0209}

\nocite{bd2}

\nocite{hmp2} \nocite{hmp3} \nocite{hmp4}

\nocite{bae1}

\end{document}